\documentclass[11pt]{article}
\usepackage{fullpage}
\usepackage{amsmath,amssymb}
\usepackage{enumerate}
\usepackage{amsmath,amssymb,amsthm}
\usepackage{graphics}
\usepackage{graphicx}
\usepackage{subfigure}
\usepackage{listings}
\usepackage[linesnumbered, ruled]{algorithm2e}

\def\uhull{{\rm uphull}}
\def\lhull{{\rm lowhull}}
\def\lema2{{\rm lemma3}}

\newtheorem{theorem}{Theorem}
\newtheorem{observation}{Observation}

\newtheorem{corollary}[theorem]{Corollary}

\newtheorem{claim}[theorem]{Claim}

\newcommand{\answerCommand}{}%
  {\renewcommand{\answerCommand}{#1\\}%
   \noindent\textbf{\answerCommand}%
   }{\\}%

  {\renewcommand{\answerCommand}{#1}%
   \noindent\textbf{\answerCommand}%
   }{\\}%

\title{Coloring geometric hypergraph defined by an arrangement of half-planes}
\author{Radoslav Fulek\thanks{Ecole Polytechnique F\'ed\'erale de Lausanne. Email:~\texttt{radoslav.fulek@epfl.ch}.
The author gratefully acknowledge support from the Swiss National
Science Foundation Grant No. 200021-125287/1.} \thanks
{This work came out of the GWOP 2009 workshop organized by Group Emo Welzl, ETH Z\"{u}rich, Switzerland}% \and
}
\date{}

\begin{document}

\maketitle

\thispagestyle{empty}

\begin{abstract}
We prove that any finite set  of half-planes can be colored by {\em two}
colors so that every point in the plane, which belongs to at least {\em
three} half-planes in the set, is covered by half-planes of both colors. This
settles a problem of Keszegh.
\end{abstract}
%\newpage

%\pagenumbering{arabic} \setcounter{page}{1}

\section{Introduction}

A {\it hypergraph} $H=(V,E)$ is a set $V$ together with a system of sets $E$ whose elements, called {\it hyperedges}, are subsets of the set $V$.
A {\it $k$-coloring} of $H$ is a mapping $\chi:V \rightarrow C$, where $|C|=k$.
We say that  an edge $e$ is {\it monochromatic} under the coloring $\chi$ if $\chi(v)$ is the same for all vertices in $e$.
A coloring $\chi$ under which no hyperedge of $H$ is monochromatic is called a {\it good} coloring.
We say that $H$ can be $k$-colored if there is a good k-coloring of $H$.
Then we define the {\it chromatic number} of $H$ to be the minimum $k$ such that $H$ can be $k$-colored.

We are concerned with specific hypergraphs obtained from half-plane arrangements. Let $\mathcal{H}$ be a finite set of half-planes in $\mathbb{R}^2$.
The set $\mathcal{H}$ defines the hypergraph $H=H(\mathcal{H})=(V,E)$ having $\mathcal{H}$ as the set of vertices, and whose hyperedges correspond to the set of points covered by at least three half-planes in $\mathcal{H}$.
More formally, for each point $p\in \mathbb{R}^2$ covered by at least three half-planes in $\mathcal{H}$, the hyperedge $e_p\in E$ is the set of half-planes $\mathcal{H}$ containing $p$.
Notice that all the points belonging to the same region in the arrangement of  lines, which define half-planes in $\mathcal{H}$, correspond to the same hyperedge (or
no hyperedge).
In \cite{Keszegh},  Keszegh  showed that the analogous hypergraph for points covered \emph{four} or more times can always be 2-colored, and asked if coverage of 3
was actually enough. We answer this affirmatively:

\begin{theorem}
\label{thm:Main}
 For any finite set of closed half-planes $\mathcal{H}$ the chromatic number of $H(\mathcal{H})$ is at most two.
 Moreover, a good  2-coloring can be computed in deterministic time $O(|V|\log |V|)$.
\end{theorem}

For the computational part of the above theorem we use the standard random access machine
theoretical model, in which every basic algebraic operation $(+,-,*,/)$ is assumed to be carried
out in a constant time.
We note that the non-algorithmic part of Theorem~\ref{thm:Main} cannot be improved: A simple example shows that
if hyperedges in  $\mathcal{H}$  correspond to the set of points covered by at least \emph{two} half-planes in $\mathcal{H}$,
its chromatic number may be at least three. In~\cite{Keszegh} it was proved that it is always at most three.

The general problem of coloring hypergraphs is well-studied and its investigation
can be traced back to the 1970s.
We note that in general it is NP-hard to
decide whether a given hypergraph is 2-colorable.
The same holds even if we restrict ourselves to 3-regular hypergraphs \cite{ll}. Hence, probably there is no
nice characterization of 2-colorable hypergraphs, if we require all hyperedges to have at least three vertices, which is our case.
Two well-known conditions for a hypergraph $H$, which are easy to check, and  which imply 2-colorability, are (1) $H$ is balanced, (2)
any union of $m$ hyperedges contains at least $m+1$ vertices (see e.g. \cite{lll}). However, neither of them can used to easily prove Theorem~\ref{thm:Main}.

Note that one can rephrase our problem in the setting of \emph{covering decomposition}.
For some recent results in the area see e.g.~\cite{Pach,Pach2}.
Thus, we can say that we want to divide $\mathcal{H}$ into two parts so that any point $p$ in the plane covered by
at least three elements of $\mathcal{H}$ is covered by a half-plane in each part.
Hence, an immediate consequence of Theorem~\ref{thm:Main} is the following.

\begin{corollary}
\label{col:Main}
 Every 3-fold covering of the plane by a finite set of closed half-planes is decomposable into two covers.
\end{corollary}

\section{Preliminaries}

From now on let $\mathcal{H}$  denote a finite set of closed half-planes in $\mathbb{R}^2$ in \emph{general position}:
no half-plane in $\mathcal{H}$ is defined by a vertical line, no two half-planes in $\mathcal{H}$ are defined
by two parallel lines, and  no three half-planes in $\mathcal{H}$ are defined by three lines intersecting in a common point.
By a standard perturbation argument one can show that general position holds without loss of generality for Theorem~\ref{thm:Main}.

We say that a half-plane in $\mathbb{R}^2$ is {\it upper} ({\it lower}),
 if it is defined as a set of points $(x,y)$, satisfying $y\leq ax+b$  (resp. $y\geq ax+b$), for some $a,b\in \mathbb{R}$.
We partition $\mathcal{H}$ into two parts  $\mathcal{H}_U$ and  $\mathcal{H}_L$ containing upper and lower half-planes, respectively.

The \emph{point-line duality} in the plane transforms the point $(a,b)\in \mathbb{R}^2$ to the line $y=ax-b$ and
the line $y=ax+b$ to the point $(a,-b)$. This duality preserves point-line incidence and above-below relationship, i.e. if a point $p$ lies above (resp. below) a line $l$,
the dual of $l$ is the point that lies below (resp. above) the line which is the dual of $p$.
The \emph{dual} of a half-plane $h$ defined by the line $y\leq ax+b$ (resp. $y\geq ax+b$) in the point-line duality is the vertical ray $r$
starting at $(a,-b)$ having downward (resp. upward) direction. This extension of the duality is natural, since a point $p\in h$, if and only if its dual line  intersects
$r$.

\begin{figure}[t]
\centering
\includegraphics[scale=0.35]{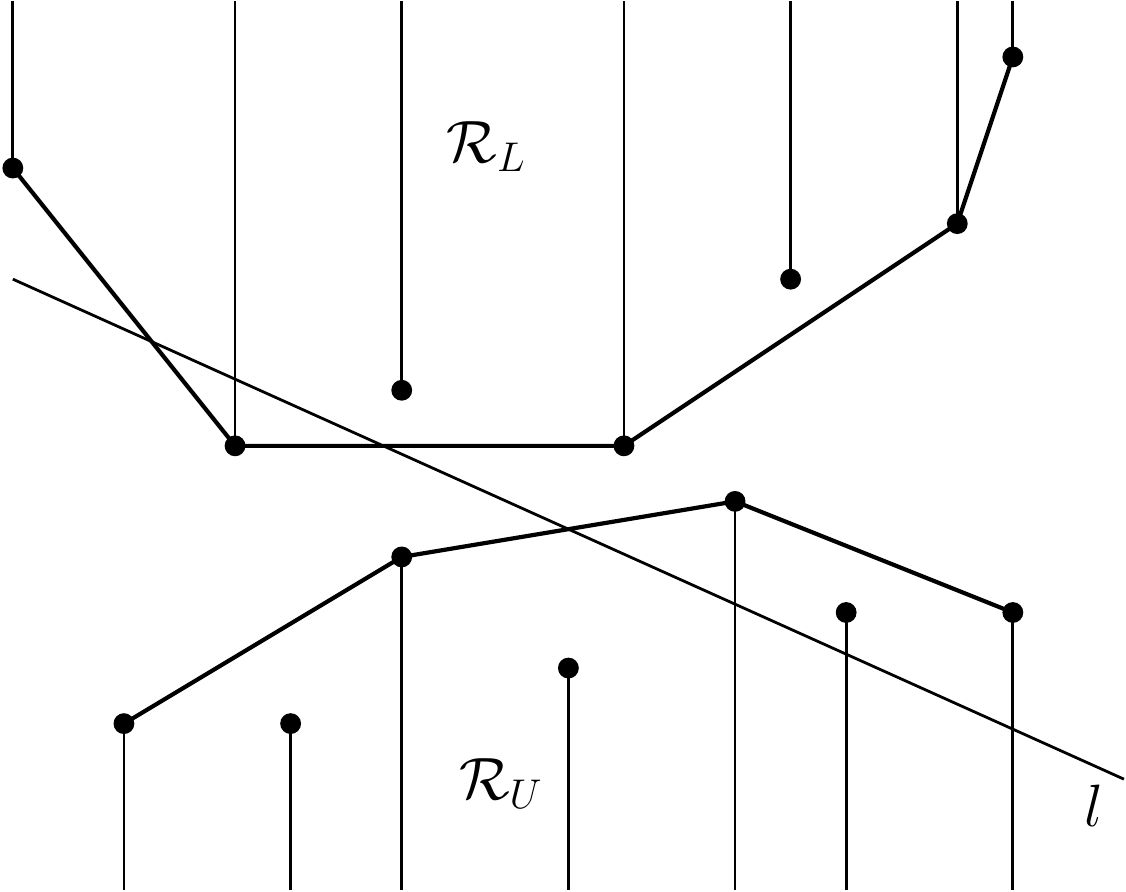}
\hspace{5mm}	
\caption{In dual settings we want to 2-color the rays so that any line $l$ intersecting at least three rays
intersects rays of both colors.}
\label{fig:first}
\end{figure}

Let $\mathcal{R}_U$  (resp.~$\mathcal{R}_L$) denote the set consisting of the rays which are duals of the half-planes in $\mathcal{H}_U$ (resp.~$\mathcal{H}_L$).
Let $\mathcal{R} = \mathcal{R}_U \cup \mathcal{R}_L$.
Using the point-line duality we can naturally recast our coloring problem so that instead of half-planes
we color the vertical rays in $\mathcal{R}$ and we require that any line $l$ intersecting at least three rays
intersects rays of both colors (see Figure~\ref{fig:first}).

 Let $\mathcal{P}_U$ (resp.~$\mathcal{P}_L$) denote the sets of starting points
of the rays in  $\mathcal{R}_U$ (resp.~$\mathcal{R}_L$). Let $\mathcal{P} = \mathcal{P}_U \cup \mathcal{P}_L$. Note that $\mathcal{P}_U$ and $\mathcal{P}_L$, respectively, could be also defined as the sets of points which are duals of
the lines defining the half-planes in $\mathcal{H}_U$ and  $\mathcal{H}_L$, respectively.
The \emph{upper} (resp. \emph{lower}) \emph{convex hull} of a set of points is the convex hull of the vertical rays directed downward (resp. upward) emanating from the points in the set.
We denote by $\mathcal{P}_U^0$ (resp. $\mathcal{P}_L^0$) the set of  vertices on the upper (resp. lower) convex hull of $\mathcal{P}_U$ (resp. $\mathcal{P}_L$).
Having defined $\mathcal{P}_U^i$ and  $\mathcal{P}_L^i$ we define $\mathcal{P}_U^{i+1}$ (resp. $\mathcal{P}_L^{i+1}$) as the set of vertices of the upper (resp. lower) convex hulls of $\mathcal{P}_U\setminus \bigcup_{j\leq i}
\mathcal{P}_U^{j}$
 and  $\mathcal{P}_L\setminus \bigcup_{j\leq i}
\mathcal{P}_L^{j}$, respectively.

Let $p,q$ denote two points in the plane. We say that $p<q$ if the $x$-coordinate of $p$ is smaller than the $x$-coordinate of $q$.

Let $\mathcal{P}'$ denote a finite set of points in $\mathbb{R}^2$.
Let  \emph{uphull}$(\mathcal{P}')$ and \emph{lowhull}$(\mathcal{P}')$, respectively, denote the
upper and lower convex hull of $\mathcal{P'}$ and  $\mathcal{P'}$, respectively.

\begin{observation}
\label{obser:CaseLemma}
If $\lhull(\mathcal{P}_L)$ and $\uhull(\mathcal{P}_U)$  intersect, then at least one of the following two sets is not empty: $\lhull(\mathcal{P}_L) \cap \mathcal{P}_U$ and  $\uhull(\mathcal{P}_U) \cap \mathcal{P}_L$.
\end{observation}

%\begin{lemma}
%\label{lemma:CaseLemma}
%Suppose that union of the elements in $\mathcal{H}$ equals to $\mathbb{R}^2$. Then
%$V(G^0)$ has a subset $S$ of size at most 2 such that for any edge $e\in E(G_0)$ of size 2,
%which has a non-empty intersection with $\mathcal{P}_U^0$ and $\mathcal{P}_L^0$, $S \cap e \not= \emptyset$.
%\end{lemma}
%\begin{proof}
%\end{proof}

Since we assume that $\mathcal{H}$ is in general position we can identify the half-planes in $\mathcal{H}$ with the points in $\mathcal{P}$. Thus,
instead of coloring half-planes in $\mathcal{H}$ (resp. rays in $\mathcal{R}$) we will  color points in $\mathcal{P}$, so that every line
intersecting at least three rays in $\mathcal{R}$ intersects rays whose corresponding points in $\mathcal{P}$ received both colors.
Our algorithm colors points incrementally. Initially all points are uncolored (depicted as squares in our pictures). One time, we assign $blue$ or $red$ colors to points
(depicted as solid black discs or empty circles, respectively).
A grey area depicts the region that does not contain any point from either $\mathcal{P}_U$ or $\mathcal{P}_L$ (depending on the situation) in its interior.

We conclude preliminaries with simple observations that serve as the main ``sub-routines'' in our algorithm for coloring the points of $\mathcal{P}$.
Intuitively, the case when $\uhull(\mathcal{P}_U)$ can be separated from $\lhull(\mathcal{P}_L)$ by a vertical line should be easier. We show
indeed that such a situation can be exploited.

\begin{figure}[h]
\centering
\includegraphics[scale=0.6]{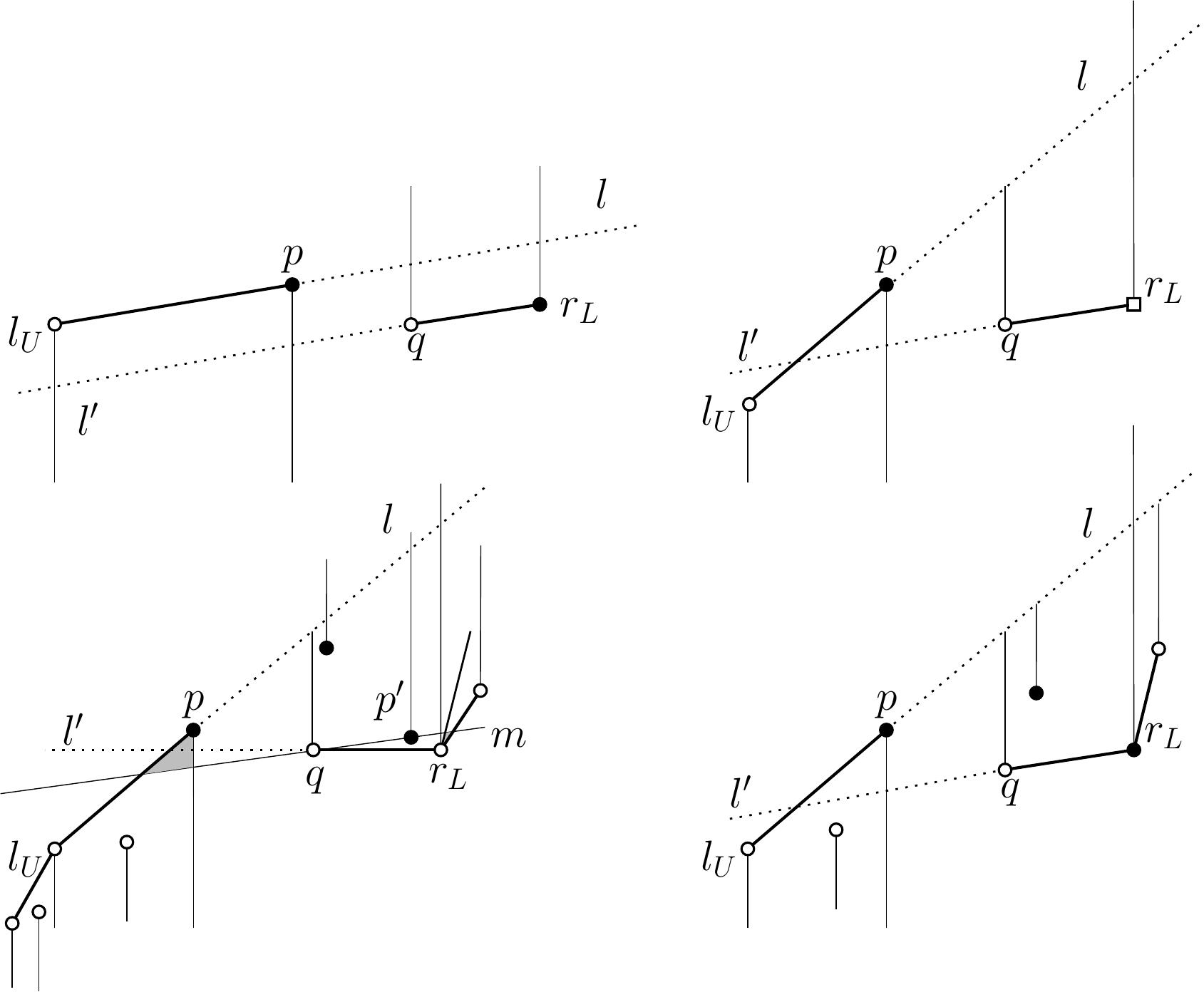}
\caption{Coloring separated hulls, red points are empty circles and blue points are solid black discs.}
\label{fig:mainLemma}
\end{figure}

Ref. to Figure~\ref{fig:mainLemma}.
Let $p\in \mathcal{P}_U^0$, and  $q\in \mathcal{P}_L^0$, $p<q$. Let $l_U$ and $r_L$ denote the points (we assume that they exist unless stated
otherwise)
 preceding and succeeding $p$ and $q$ on
their respective hulls.
Suppose that $\mathcal{P}_U$ does not contain any point to the right of $p$, and
$\mathcal{P}_L$ does not contain any point to the left of $q$
Moreover, we assume that the line $l$ through $l_U$ and $p$ passes above $q$, and
 the line $l'$ through $q$ and $r_L$ passes below
$p$.

Under the assumption of the previous paragraph.

%Let $r_U$ denote the point succeeding $p$ on the upper hull (if it exists). Let $l_L$ denote the point succeeding $q$ on the lower hull (if it exists).

\begin{observation}
\label{lemma:MainLemma1}
%Suppose we have a good 2-coloring $\chi$ of all points in $\mathcal{P}$, except all $r\in \mathcal{P}_U$, $r<p$ and
%$r\in \mathcal{P}_L$, $r>q$, so that $\chi(p)\not=\chi(q)$, and if the interior of the triangle $l_Upr_U$ $(l_Lqr_L)$   does not contain a point from
%$\mathcal{P}_U$ ($\mathcal{P}_L$), we have $\chi(p)\not=\chi(r_U)$ $(\chi(q)\not=\chi(l_L))$.
If neither $l$ intersects the segment $qr_L$, nor $l'$ intersects the segment $l_Up$, then
the 2-coloring of $\mathcal{P}$, which colors $p$ with $blue$, $q$ with $red$,
all the vertices in $\mathcal{P}_{U}$ to the left of $p$ by $red$, and all the vertices in $\mathcal{P}_{L}$
to the right of $q$ by $blue$, is a good coloring.
\end{observation}

See the upper left part of Figure~\ref{fig:mainLemma} for the proof of Observation~\ref{lemma:MainLemma1}.

Under the same assumption as in Observation~\ref{lemma:MainLemma1}.

\begin{observation}
\label{lemma:MainLemma2}
%Suppose we have a good 2-coloring $\chi$ of all points in $\mathcal{P}$, except all $r\in \mathcal{P}_U$, $r<p$ and
%$r\in \mathcal{P}_L$, $r>q$, so that $\chi(p)\not=\chi(q)$, and if the interior of the triangle $l_Upr_U$ $(l_Lqr_L)$   does not contain a point from
%$\mathcal{P}_U$ ($\mathcal{P}_L$), we have $\chi(p)\not=\chi(r_U)$ $(\chi(q)\not=\chi(l_L))$.
If $l'$ intersects the segment $l_Up$ or $l_U$ does not exist,
there exists a good 2-coloring of $\mathcal{P}$, which colors $p$ with $blue$, $q$ with $red$,
all the vertices in $\mathcal{P}_{U}$  to the left  of $p$  by $red$, all the vertices in $\mathcal{P}_{L}$
 between $q$ and $r_L$  by $blue$, and all the vertices to the right of $r_L$  by $red$. Analogously, if
 $l$ intersects $qr_L$ or $r_L$ does not exist, there is a good 2-coloring of $\mathcal{P}$.
\end{observation}

\begin{proof}
The situation is depicted in the upper right part of Figure~\ref{fig:mainLemma}.
Observe that we need only to decide the color of $r_L$.

Let $m$ denote the tangent through $q$ to $\lhull(\mathcal{P}_L^1)$, if $\mathcal{P}_L^1\not=\emptyset$. If $m$ passes through a point $p'$ of $\mathcal{P}_L^1$ between $q$ and $r_L$, $m$ is below all the points of $\mathcal{P}_L$ except $r_L,p'$ and $q$, and $m$ is above all the points in $\mathcal{P}_U$ except $p$, we
color $r_L$ by $red$ (see the lower left part of Figure~\ref{fig:mainLemma}). Otherwise, we color $r_L$ by $blue$ (see the lower right part of Figure~\ref{fig:mainLemma}).

It is straightforward to check that our 2-coloring is good in both cases. Indeed, there cannot be a $red$ monochromatic edge, and
a $blue$ monochromatic edge has to contain $r_L$ and $p$, which yields the tangent $m$  with the above properties and that in turn implies the $red$ color for $r_L$.
\end{proof}

%\newpage

\section{Proof of Theorem \ref{thm:Main}}

First, we deal with the case when there exists a point $p\in \mathbb{R}^2$ that is not covered by any half-plane in $\mathcal{H}$.

\begin{claim}
\label{claim:covered}
We can assume that every point in $\mathbb{R}^2$ is contained in at least one half-plane of $\mathcal{H}$.
\end{claim}

\begin{proof}
Assume that the origin $o=(0,0)$ is not covered by any half-plane in $\mathcal{H}$.
The \emph{point-line polar duality} in the plane transforms the line $ax+by=1$, $(a,b)\not= (0,0)$ to the point $(a,b)$ and vice versa.

We reduce our problem using the point-line polar duality to a problem of coloring a hypergraph $H'=(V',E')$ defined as follows.
The set of vertices of $H'$ is a finite set of points in the plane and a hyperedge in $E'$ is the intersection of a closed half-plane with $V'$
of size at least three. In \cite{Keszegh} it was shown that $H'$ can be always two-colored by an algorithm with the running time of $O(|V'|\log |V'|)$.

 We use the
polar duality on the lines defining the half-planes in $\mathcal{H}$ thereby obtaining a set of points $\mathcal{P}_o$. Let $\mathcal{L}_o$ denote
the set of line segments $op$, where $p\in \mathcal{P}_o$.
Now, it is enough to two-color the line segments in $\mathcal{L}_o$ so that any line intersecting at least three line segments in $\mathcal{L}_o$ intersects
line segments of both colors.
 We use the algorithm from \cite{Keszegh} to two-color the points in $\mathcal{P}_o$.
A good coloring of the line segments in $\mathcal{L}_o$ is obtained by assigning to every line segment the color of its endpoint in $\mathcal{P}_o$.
\end{proof}

Thus, by Claim~\ref{claim:covered} we can assume that the whole plane is covered by the half-planes in $\mathcal{H}$.

%Hence, we can use  Lemma  \ref{lemma:CaseLemma}, thereby obtaining a set $S$, $|S|\leq 3$.
Note that the assumption about covering the plane by the half-planes in $\mathcal{H}$ translates in the dual setting
to the assumption that $\uhull(\mathcal{P}_U)$ and $\lhull(\mathcal{P}_L)$ intersect.
Thus, by using Observation \ref{obser:CaseLemma} we obtain a point $p\in \mathcal{P}_U^0$ (w.l.o.g.) contained in $\lhull(\mathcal{P}_L)$.
Hence, we have two points $l_L,r_L\in \mathcal{P}_U^0$, $l_L<p<r_L$, such that there is no point $q\in \mathcal{P}_L^0$, for which $l_L<q<r_L$.
By left-right symmetry we can assume that $p$ is not the leftmost point in $\mathcal{P}_U$ unless $|\mathcal{P}_U|=1$.\\

Let us assume that $|\mathcal{P}_U|>1$ (the case when $|\mathcal{P}_U|=1$ is discussed later).  Let $l_U\in \mathcal{P}_U$ denote the point immediately to the left of $p$ on $\uhull(\mathcal{P}_U)$,
and let $h$ denote the line through $p$ and $l_U$.
Let $r_U\in \mathcal{P}_U$ denote the point immediately to the right of $p$ on the upper hull (if it exists).
 Let $v$ denote the vertical line through $p$.
% The lines $v$ and $h$ divide the plane into 4 regions (see Figure \ref{fig:caseb}).
Depending on the position of $l_L$ and $r_L$ above or below $h$, on the existence of an intersection between the segments $l_Lr_L$ and $l_Up$, and on whether $l_L<l_U$ holds, we will distinguish the following 4 cases (a)-(d).

In each case we define a good 2-coloring $\chi$ of $H$.
%A triangles and discs positioned on a vertical ray starting at some vertex of $\mathcal{P}$ below its starting point mean that
%its  corresponding vertex  might be colored with some other color depending on a particular situation.
%A triangles and discs not positioned on a vertical ray depicts a coloring of  the points in a particular region.

\begin{enumerate}[a)]

\item
\label{case:b}
\begin{figure}[t]
\centering
\includegraphics[scale=0.6]{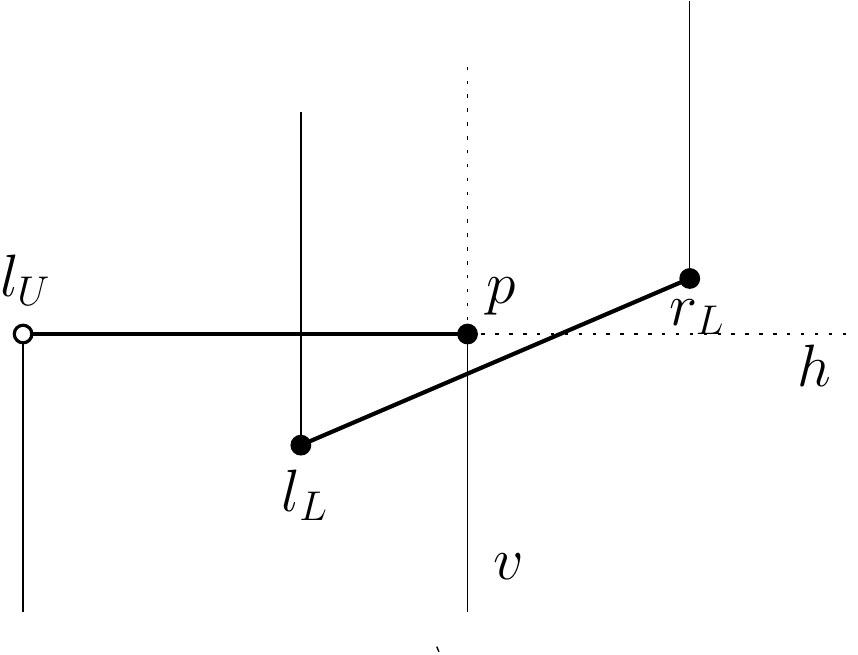}
\caption{Case (\ref{case:b})}
\label{fig:caseb}
\end{figure}

In this case we have: $r_L$ is above $h$, which implies that  $l_Up$ and $l_Lr_L$ do not intersect each other (see Figure \ref{fig:caseb}).

We color the points as follows: $\chi(p)=\chi(r_L)=\chi(l_L)=blue$, and the remaining points by $red$.

The coloring is good as any non-vertical line intersects a ray corresponding to $p,r_L$ or $l_L$, and no
line can intersect all the rays corresponding to $p,r_L$ and $l_L$ without intersecting the ray corresponding to $l_U$.

\item
\label{case:c}

\begin{figure}[t]
\centering
\includegraphics[scale=0.6]{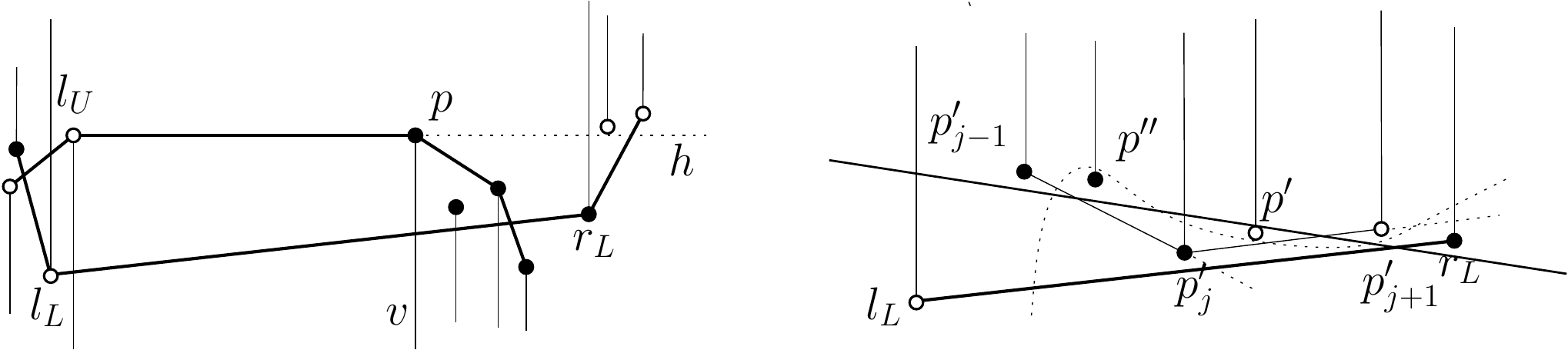}
\caption{Case (\ref{case:c})}
\label{fig:casec}
\end{figure}

In this case we have:  $r_L$ is below $h$, $l_Up$ and $l_Lr_L$ do not intersect each other, and
$l_L<l_U$ (see Figure~\ref{fig:casec} left).

We color the points as follows: $\chi(p)=\chi(r_L)=blue$ and $\chi(l_U)=\chi(l_L)=red$.
We color the points $q$;
$q\in \mathcal{P}_U$, $q<l_U$ or $q\in \mathcal{P}_L$, $q>r_L$, by $red$ and the remaining points $q$;
$q\in \mathcal{P}_U$, $q>p$ or  $q\in \mathcal{P}_L$, $q<l_L$, by $blue$.
The points in $\mathcal{P}_U$ between $l_U$ and $p$ can be colored arbitrarily.
Now, we describe the coloring of the rest of $\mathcal{P}_L$.

Let $p_1',p_2',\ldots ,p_i'\in \mathcal{P}_L^1$, $p_1'<p_2'<\ldots<p_i'$, denote the points between $l_L$ and $r_L$ (if they exist) in $\mathcal{P}_L^1$.
Let $p_0'<p_1'$ (resp. $p_{i+1}'>p_i'$) denote the point in  $\mathcal{P}_L^{1}$ immediately preceding $p_1'$ (resp. following $p_i'$).
If $p_0'$ (resp. $p_{i+1}'$) does not exist let $p_0'$ (resp. $p_{i+1}'$) denote a point below $p_1'$ slightly to the left (resp. below $p_i'$ slightly to the right).
Let $p_j'$, $j>0$, denote the point with minimal $j$ so that the line through $p_{j+1}'$ and $p_j'$ passes above $r_L$.
All the points in $\mathcal{P}_L$  between $l_L$ and $p_j'$ are colored with $blue$ and between $p_j'$ and $r_L$ with $red$.
The color of $p_j'$ is defined as follows.
If $p_j'$ forms a hyperedge only with $l_L$, and a point in $\mathcal{P}_L$  between $p_j'$ and $r_L$, we color $p_j'$ with $blue$ (see Figure~\ref{fig:casec} right). Otherwise, we color $p_j'$ with $red$. The last condition can be also expressed as follows: If there exists a line through $r_L$ and through a point in $\mathcal{P}_L$ between
$p_j'$ and $r_L$ passing above $l_L$ and $p_j'$, and passing below all the other points in $\mathcal{P}_L$ between $l_L$ and $r_L$,
we color $p_j'$ with $blue$. Otherwise, we color $p_j'$ with $red$.

%We color $p_1'$ by $blue$ and $p_3'$ by $red$. If there is a line $l$ such that $l(l_L,p_2',p_3')$, $\neg l(p_1')$ and $\neg l(r_L)$,
%we color $p_2'$ with $blue$ (see Figure \ref{fig:casec} right), otherwise we color $p_2'$ with $red$.
%The rest of the points to the left of  $p_2'$ is colored by $blue$ and to the right of $p_2'$ by $red$.

In what follows we check that the 2-coloring we defined is good:
Any line $l$ witnessing a $blue$ monochromatic edge $e$ has to pass below $l_L$ and above $l_U$. Moreover, $l$ has to pass below all the vertices of $\mathcal{P}_L^0$ except $r_L$, and below all the vertices of $\mathcal{P}_L$ to the right of $r_L$. Hence, all the points
in $\mathcal{P}$ participating in $e$ are  points from $\mathcal{P}_L$ between $l_L$ and $r_L$. By left-right red-blue symmetry, the same holds for $red$ monochromatic edges. Hence, it is enough to show that no monochromatic hyperedge is formed by points in $\mathcal{P}_L$ between $l_L$ and $r_L$, including $l_L$ and $r_L$.

Observe that it cannot happen that $p_j'$ forms both (1) a hyperedge only with $l_L$ and a point in $\mathcal{P}_L$  between $p_j'$ and $r_L$, (2) a hyperedge only with $r_L$ and a point in $\mathcal{P}_L$  between $l_L$ and $p_j'$ .
Indeed, otherwise we find two different lines intersecting in more than one point (see Figure~\ref{fig:casec} right).
Similarly, it cannot happen that a hyperedge is formed only by $r_L$, and points in $\mathcal{P}_L$ between $l_L$ and $p_j'$, since
that would violate the minimality of $j$. Finally, it also cannot happen that a hyperedge is formed only by $l_L$ and points in $\mathcal{P}_L$ between $p_j'$ and $r_L$.

\item
\label{case:a}
In this case we have:  $r_L$ is below $h$, and $l_Up$ and $l_Lr_L$ intersect.

 Let $\chi(p)=blue$ and $\chi(l_L)=red$ (the color of $l_L$ might be changed in some of the following subcases).
Let $r_L'$ (resp. $l_L'$) denote the point following $r_L$  (resp. preceding $l_L$) on the lower hull (if it exists).

First, we assume that either $r_L$ is above the line $pr_U$, or $r_U$ does not exist.
Observe that in this case we can also assume that $r_U$ does not belong to $\lhull(\mathcal{P}_L)$. Indeed, otherwise $r_U$ can play the role of the point $p$ and
we easily reduce our situation to case (\ref{case:b}).
In what follows we distinguish several subcases:
\begin{enumerate}[(c1)]
%, possibly after reversing the
%$x$ or $y$-axis,
\item
The line $l$ through $r_L$ and $r_U$ (resp. $l_L$ and $l_U$) passes  below all the points in $\mathcal{P}_L$ except $l_L$ and $r_L$, and  above all the points in $\mathcal{P}_U$ except $r_U$ (resp. $l_U$) (see Figure~\ref{fig:casea1}).

We put $\chi(l_U)=red$ and $\chi(r_U)=\chi(r_L)=blue$ (resp.  $\chi(l_U)=\chi(l_L)=blue$ and $\chi(r_U)=\chi(r_L)=red$).
We color the rest of the points by $red$ (resp. $blue$).

The defined coloring is good, as a non-vertical line that does not intersect any of the rays corresponding to $p,r_L$ or $r_U$ (resp. $p,l_L$ or $l_U$)
cannot intersect any ray except the one corresponding to $l_L$ (resp. $r_L$).
Moreover, a line cannot intersect all the rays corresponding to $p,r_L$ and $r_U$ (resp. $p,l_L$ and $l_U$).

\label{case:a1}

\item

\begin{figure}[t]
\centering
\subfigure[]{\label{fig:casea1}
		\includegraphics[scale=0.4]{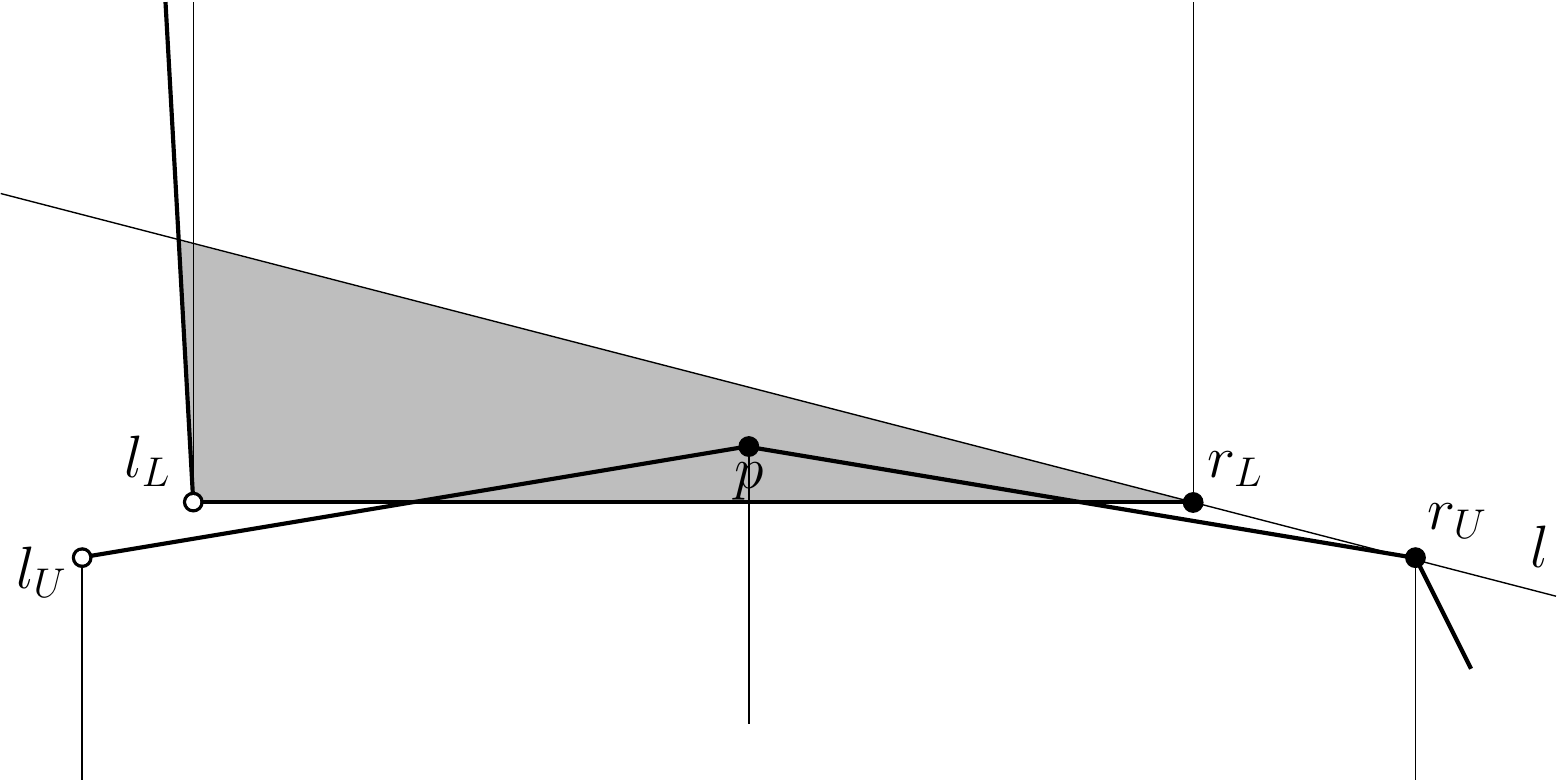}
	    \hspace{5mm}
	}
  \subfigure[]{  \label{fig:casea3}
		
		\includegraphics[scale=0.4]{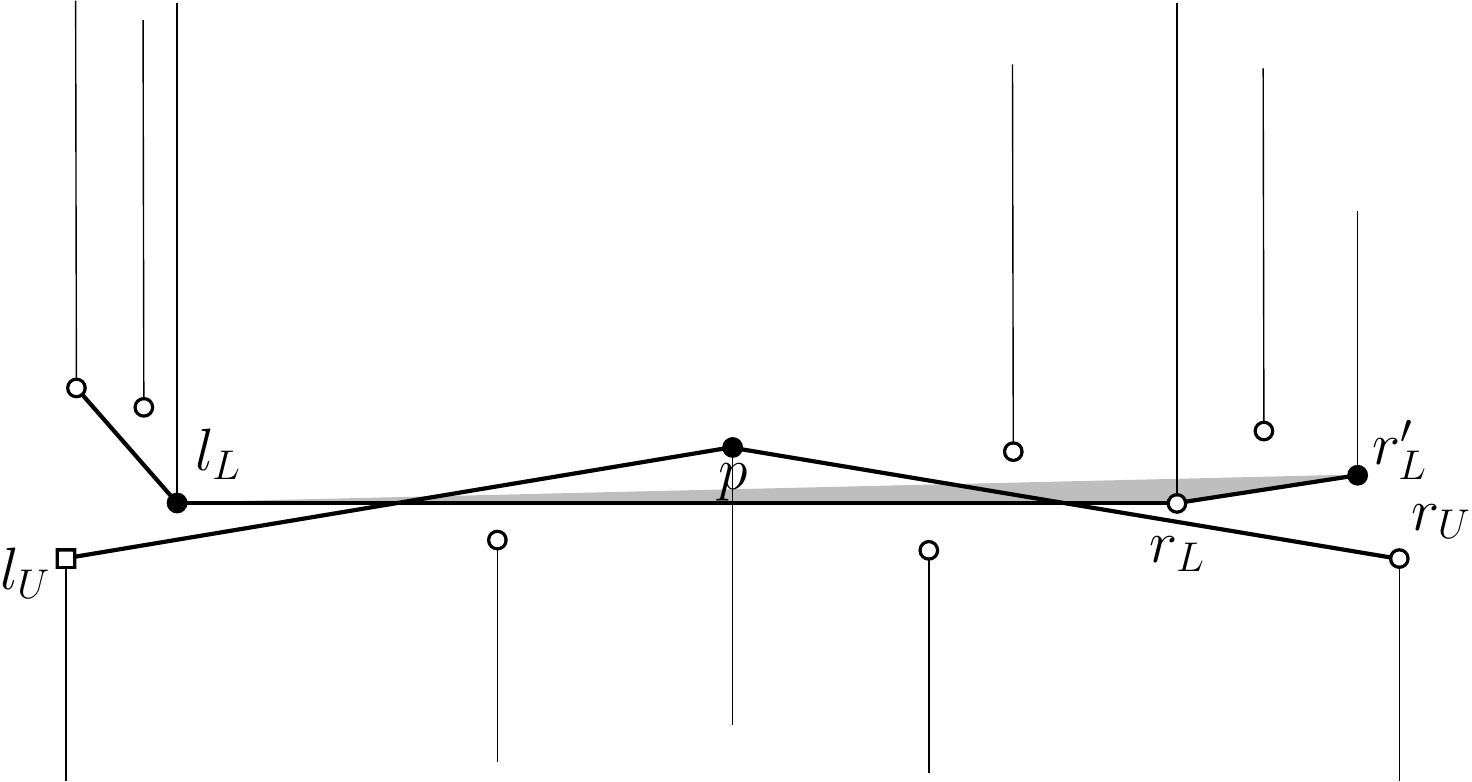}
  }
	\caption{(a)  Case (\ref{case:a1}), (b) Case (\ref{case:a3})}
\end{figure}

% \begin{figure}[t]
% \centering
% \includegraphics[scale=0.6]{FinalA0}
% \caption{Case (\ref{case:a1})}
% \label{fig:casea1}
% \end{figure}
\begin{figure}[t]
\centering
\subfigure[]{	\label{fig:picablesk}	\includegraphics[scale=0.4]{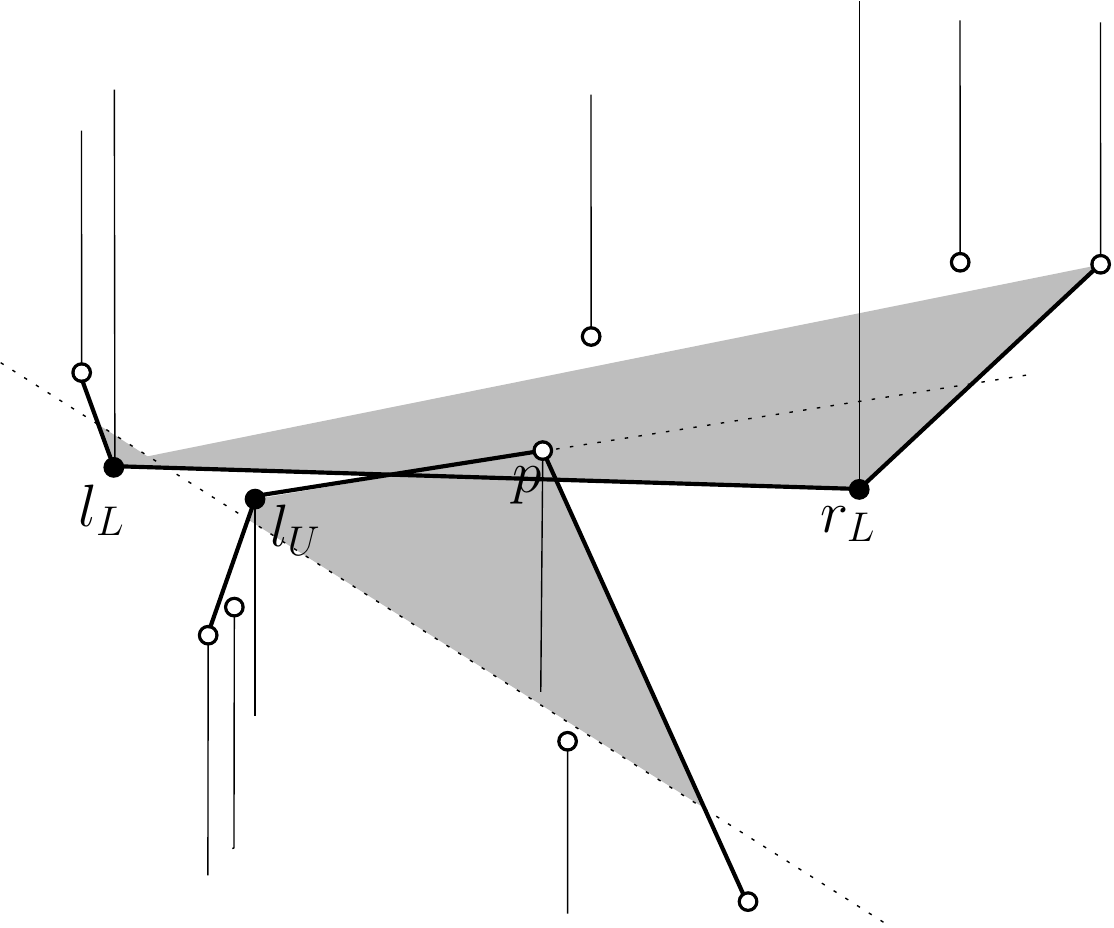}}
\hspace{5mm}
\subfigure[]{	\label{fig:picablesk2}			\includegraphics[scale=0.4]{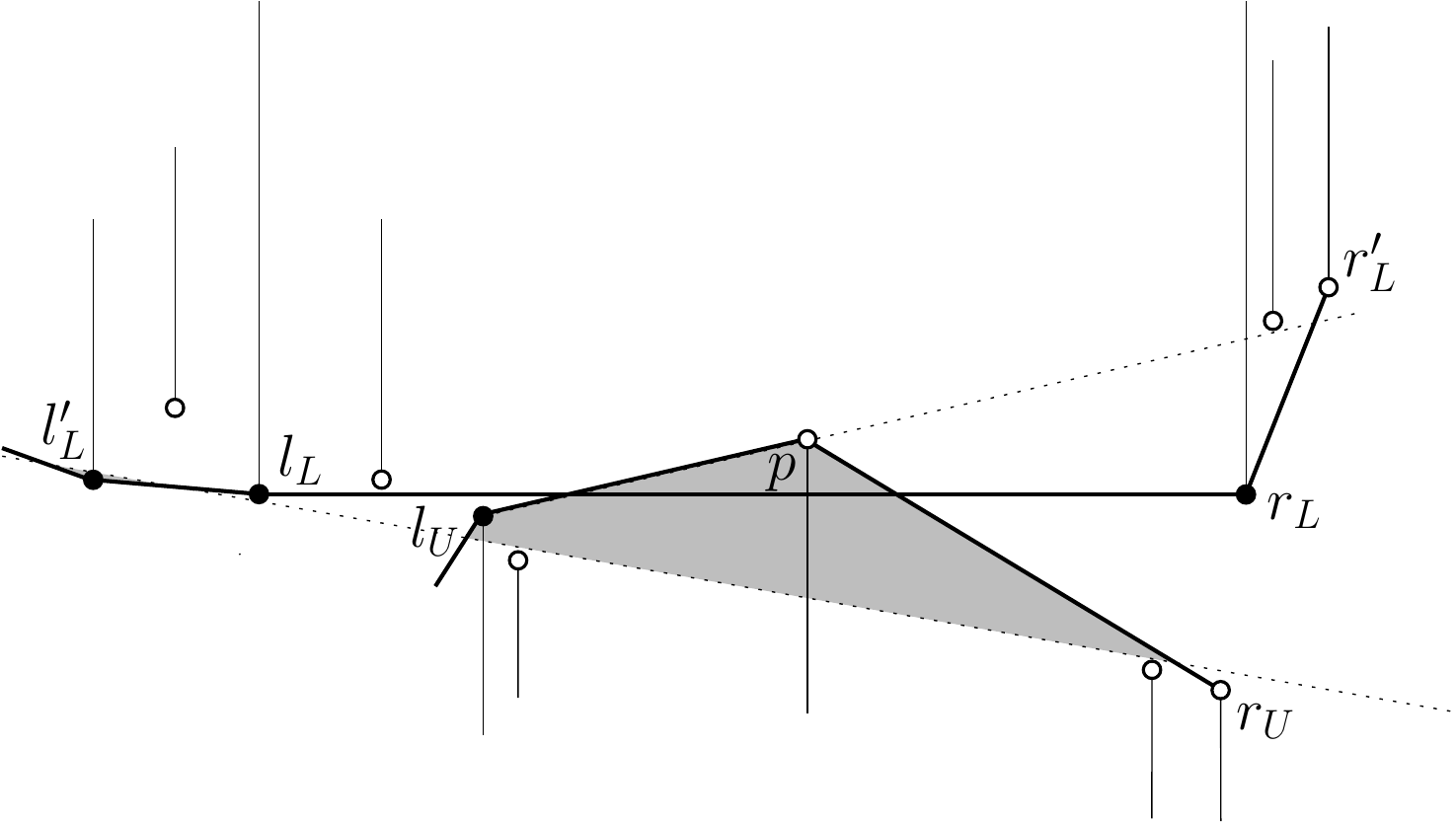}}

\caption{(a) a hyperedge formed by $l_L,l_U$ and $p$, (b) a hyperedge formed by $l_L',l_U$ and $p$}

\end{figure}

\begin{figure}[t]
\centering
\subfigure[]{\label{fig:casea2}
		\includegraphics[scale=0.4]{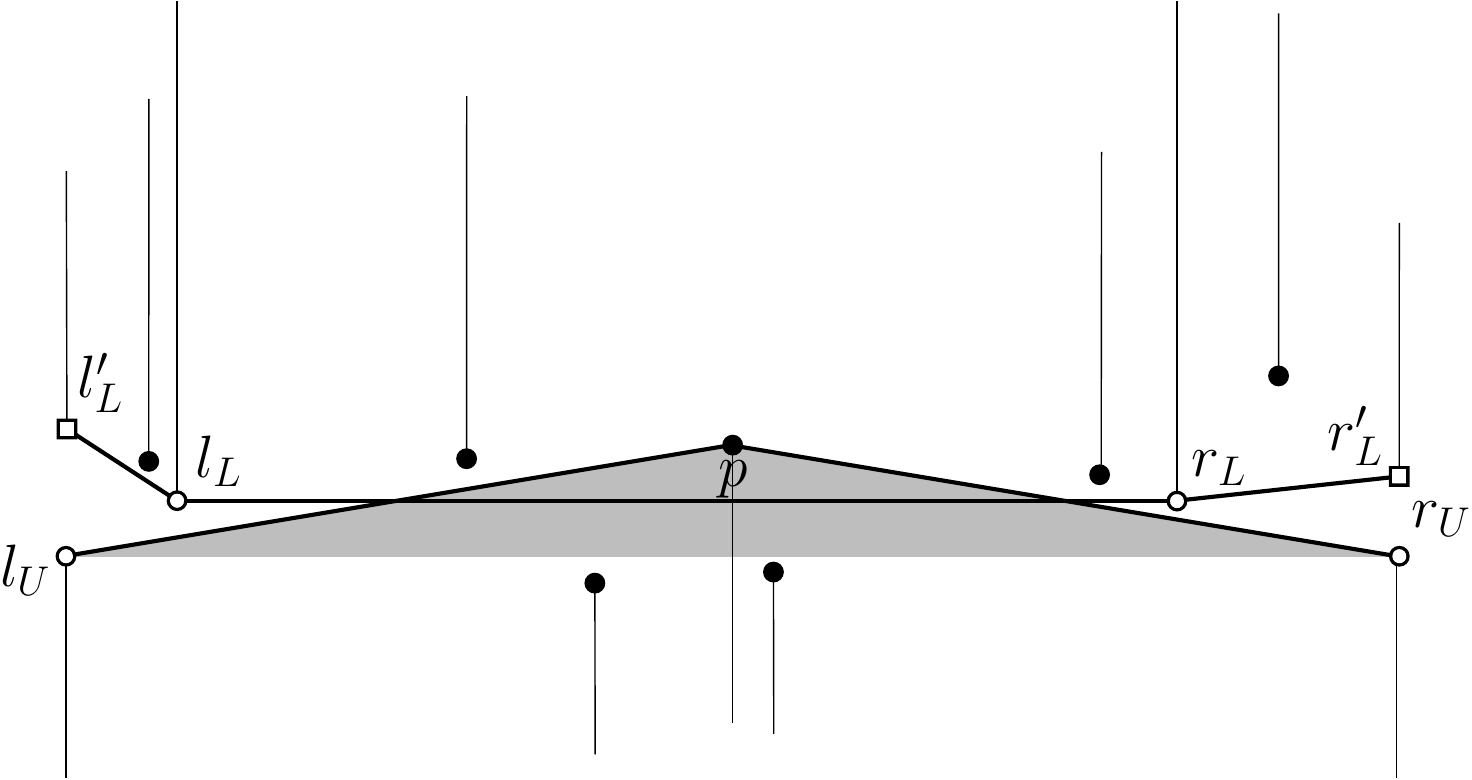}
	    \hspace{5mm}
	}
  \subfigure[]{\label{fig:casea4}
		\includegraphics[scale=0.4]{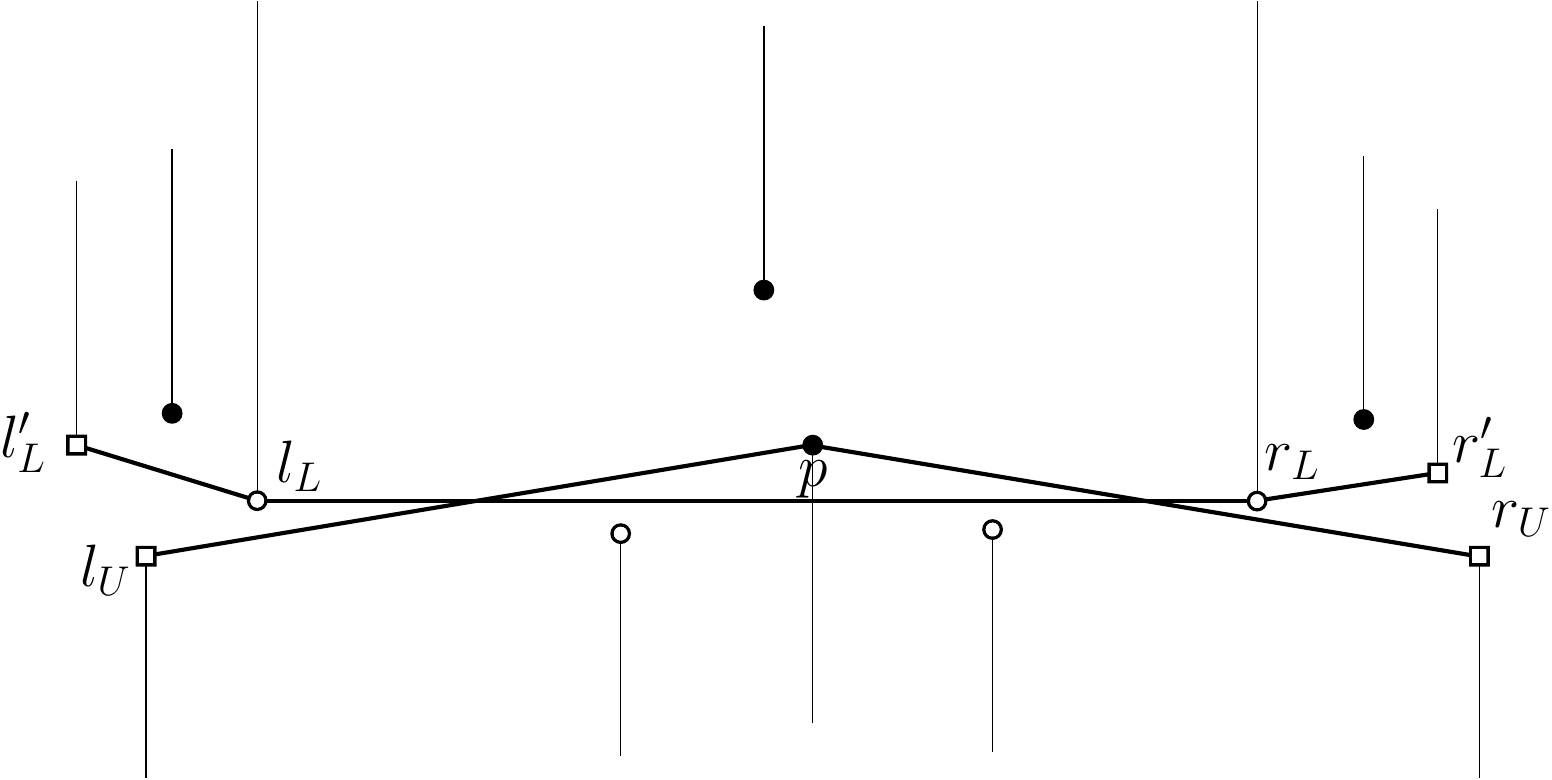}
  }
	\caption{(a)  Case (\ref{case:a2}), (b) Case (\ref{case:a4}) }
\end{figure}

The triangle $l_Lr_Lr_L'$ (resp. $r_Ll_Ll_L'$) does not contain any point from $\mathcal{P}_L$ in its interior (see Figure \ref{fig:casea3}).
We assume that $r_L'$ (resp. $l_L'$) exists.
We handle only the case concerning the triangle $l_Lr_Lr_L'$, since the situation when $r_Ll_Ll_L'$ does not contain any point from $\mathcal{P}_L$ in its interior can be handled by symmetry.

We put $\chi(l_L)=\chi(r_L')=blue$, and $\chi(r_U)=\chi(r_L)=red$.
We delete the points in $\mathcal{P}_{L}$ between $l_L$ and $r_L'$, except $r_L$, since
they can be always colored with  $red$.
We color the points $r\in\mathcal{P}_U$, $r>p$, and
 $r\in\mathcal{P}_L$, $r<l_L$, with $red$.
% We color the points in $\mathcal{P}_L$ between $l_L$ and $r_L'$ with $red$.
 %If the line through $r_L$ and $r_L'$ intersects the segment $pl_U$ we color the rest of the points by $red$.
We apply either Observation~\ref{lemma:MainLemma1} or~\ref{lemma:MainLemma2} (depending on where the line through $p$ and $l_U$
meets the line through $r_L$ and $r_L'$)
with $p$ as $p$ and $r_L$ as $q$ in order to color the rest of the points.
%, without taking the vertices in $\mathcal{P}_L$ between $r_L$ and $r_L'$ into consideration.
Note that $r'_L$ was not recolored by the observation, as the points in $\mathcal{P}_L$ between $r_L$ and $r_L'$ were
deleted. Thus, by the proof of Observation~\ref{lemma:MainLemma2} $r_L'$ is always colored $blue$.

The coloring we define in this case might not yet be good, as $l_L,l_U$ and $p$ might form a monochromatic hyperedge
(see Figure~\ref{fig:picablesk}). This is equivalent to the situation when a tangent from $l_U$ to the lower hull of $\mathcal{P}_L$ passes through $l_L$ and above all the points in $\mathcal{P}_U$ except $p$ and $l_U$. However, in this case we can color everything with $red$, except $l_L,l_U$
and $r_L$.

Let us check that our coloring is good.
A $red$ monochromatic hyperedge is rule out easily by either Observation~\ref{lemma:MainLemma1} or~\ref{lemma:MainLemma2}.
A $blue$ monochromatic hyperedge $e$ has to contain points that were not colored by the application of either  Observation~\ref{lemma:MainLemma1} or~\ref{lemma:MainLemma2}. The only such $blue$ point is $l_L$. Hence, $e$ also contains $p$ and $l_U$ (by Observation~\ref{lemma:MainLemma2}). Now, we have everything colored by $red$ except $l_L,l_U$ and $r_L$.

In this situation a $blue$ monochromatic hyperedge is easily ruled out as no line can intersect three rays corresponding to
$l_L,l_U$ and $r_L$. On the other hand the line that misses all these three rays can hit only the ray corresponding to $p$.

\label{case:a3}
\item
The triangle $l_Upr_U$ does not contain any point from $\mathcal{P}_U$ in its interior (see Figure \ref{fig:casea2}), and
none of the above happens.

We put $\chi(l_U)=\chi(r_U)=\chi(r_L)=red$. We delete all the points in $\mathcal{P}_{U}$ between $l_U$ and $r_U$, except $p$,
since they can be always colored with $blue$. We color all the points in $\mathcal{P}_{L}$ between $l_L$ and $r_L$ with $blue$.
%If the line $r_Lr_L'$ intersects the segment $l_Up$,
We apply either Observation~\ref{lemma:MainLemma1} or~\ref{lemma:MainLemma2} with $p$ as $p$ and $r_L$ as $q$
in order to color the points in $\mathcal{P}_U$ to the left of $p$ and  in $\mathcal{P}_L$ to the right of $r_L$.
Analogously, we apply  either Observation~\ref{lemma:MainLemma1} or~\ref{lemma:MainLemma2} (with the orientation of the $x$-axis reversed) with $p$ as $p$ and $l_L$ as $q$
in order to color the points in $\mathcal{P}_U$ to the right of $p$ and in $\mathcal{P}_L$ to the left of $l_L$.
Neither $l_U$ nor $r_U$ is recolored, by the coloring constructed in the proof of either Observation~\ref{lemma:MainLemma1} and~\ref{lemma:MainLemma2}.

It is easy to check that the 2-coloring we defined is good, using the fact that we excluded  the previous cases (\ref{case:a1}), (\ref{case:a3}).
Indeed, $r_U,r_L$ and $l_L$ (resp. $l_U, l_L$ and $r_L$) cannot form a hyperedge, by excluding case (\ref{case:a1}).
By excluding case (\ref{case:a3}) the triangle $l_Lr_Ll_L'$ (resp. $l_Lr_Lr_L'$) contains a point from $\mathcal{P}_L$, which must be $blue$.
Thus, $r_L,l_L$ and $l_L'$  (resp. $r_L,l_L$ and $r_L'$) also cannot form a monochromatic hyperedge.
%Nothing else "bad" can happen by the coloring constructed in the proof of Lemma \ref{lemma:MainLemma}.

\label{case:a2}
% \begin{figure}[t]
% \centering
% \includegraphics[scale=0.5]{FinalA-LPR-U-E}
% \caption{Case (\ref{case:a2})}
% \label{fig:casea2}
% \end{figure}

%By the proof of Lemma \ref{lemma:MainLemma}  $r_L'$ receives in this case  $blue$ as well.

%Note that if both $l_U$ and $r_U$ receive $blue$, everything is fine since the triangle $l_Ur_Up$ contains an element from
%$\mathcal{P}_U$, which is in any case colored with $red$ color.

% \begin{figure}[t]
% \centering
% \includegraphics[scale=0.5]{FinalA-LRR-L-E}
% \caption{Case (\ref{case:a3})}
% \label{fig:casea3}
% \end{figure}

\item

None of the previous cases occurs.

We put $\chi(r_L)=red$ (see Figure \ref{fig:casea4}).
We apply either Observation~\ref{lemma:MainLemma1} or~\ref{lemma:MainLemma2} twice: first with $p$ as $p$ and $r_L$ as $q$;
then with the orientation of the $x$-axis reversed, with $p$ as $p$ and $l_L$ as $q$.
Finally, we color all the vertices in $\mathcal{P}_{L}$ between $l_L$ and $r_L$ with $blue$.

Note that if both $l_U$ and $r_U$ receive $blue$, $l_U,p,$ and $r_U$ cannot form the monochromatic hyperedge,
 since the triangle $l_Ur_Up$ contains an element from
$\mathcal{P}_U$, which is in this case colored with $red$.
 Similarly, we can handle the situation if $r_L'$ or $l_L'$ receives $red$.
On the other hand, it can still happen that either $l_L',l_U,p$ or  $r_L',r_U,p$ forms a monochromatic blue edge analogously to
case (\ref{case:a3}); due to symmetry we treat only the first case.  This is equivalent to the following disjunctive condition:
\emph{either} a tangent through $l_U$ to the lower convex hull of $\mathcal{P}_L$ passes through $l_L'$, and above all the points in $\mathcal{P}_U$ except $l_U$ and $p$
\emph{or} the line through $l_L'$ and $l_L$ passes above all the points in $\mathcal{P}_U$ except $l_U$ and $p$.

However, by the coloring constructed in the proof of Observation~\ref{lemma:MainLemma2}, if that is the case (see Figure~\ref{fig:picablesk2}) the line through
$l_L'$ and $l_L$ intersects the segment $pr_U$. Thus, we can color
everything by $red$, except $l_L', l_L, l_U$ and $r_L$.
The coloring is still good as $l_L'l_Lr_L$ contains a point from $\mathcal{P}_L$ in its interior.
 %Nothing else bad can happen by the coloring constructed in the proof of Lemma \ref{lemma:MainLemma}.

\label{case:a4}
% \begin{figure}
% \centering
% \includegraphics[scale=0.5]{FinalA-NE}
% \caption{Case (\ref{case:a4})}
% \label{fig:casea4}
% \end{figure}
\end{enumerate}

\begin{figure}[t]
\centering
\label{fig:casea}
		\includegraphics[scale=0.4]{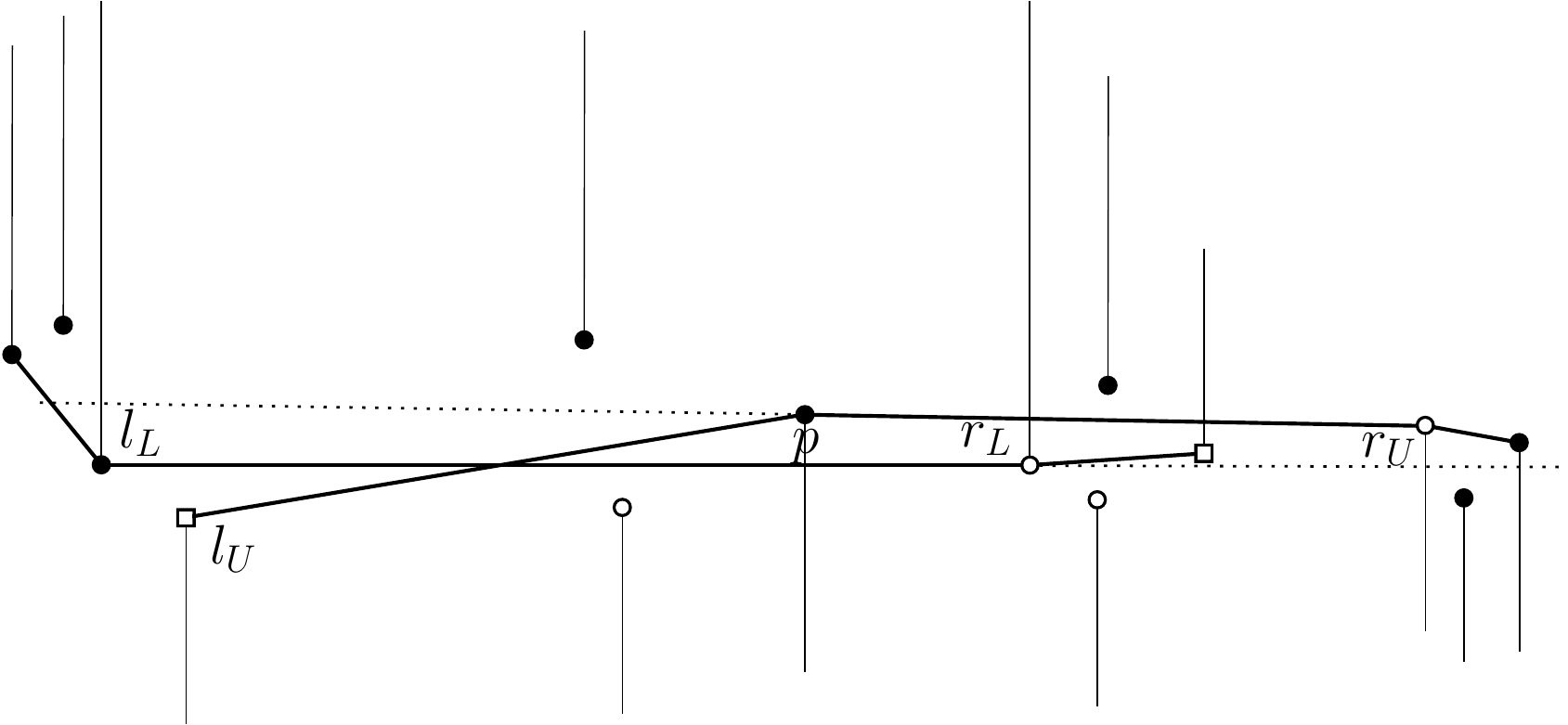}

	\caption{$r_L$ is below the line $pr_U$}
\end{figure}

\begin{figure}[t]
\centering
\subfigure[]{\label{fig:cased0}
		\includegraphics[scale=0.4]{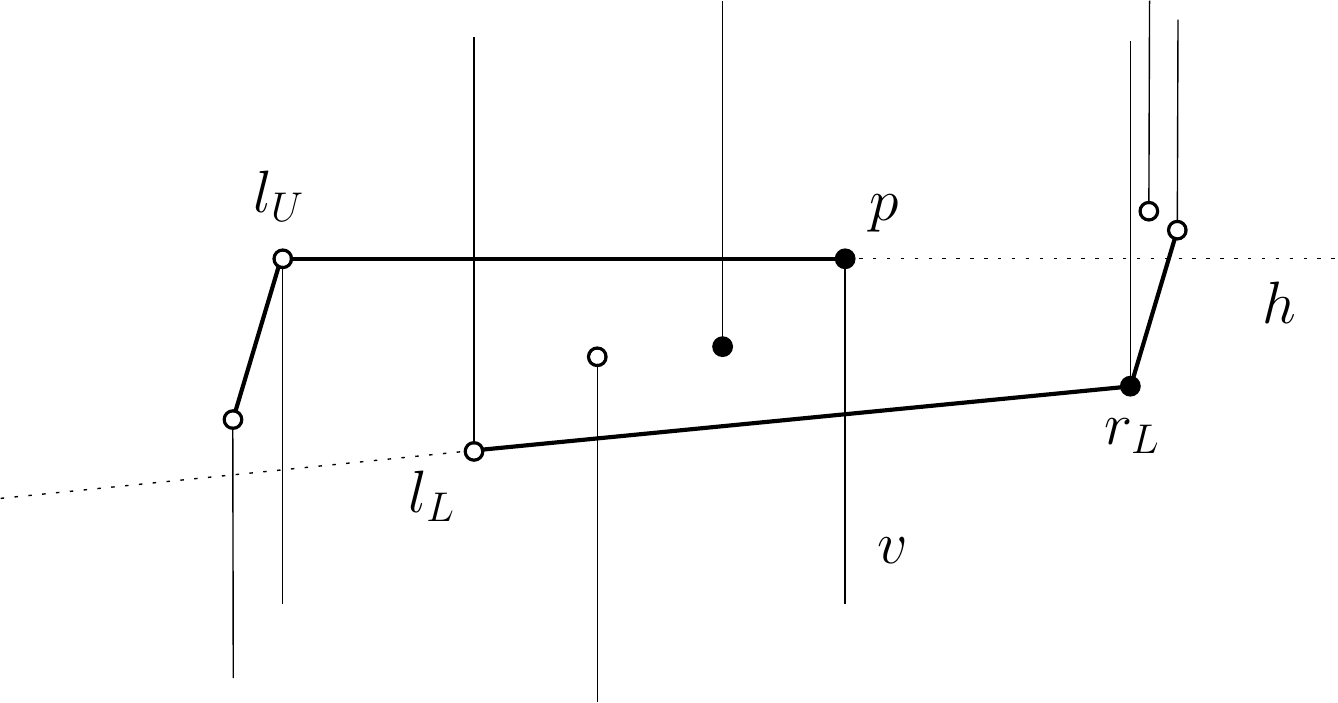}
	    \hspace{5mm}
	}
  \subfigure[]{\label{fig:cased}
		\includegraphics[scale=0.4]{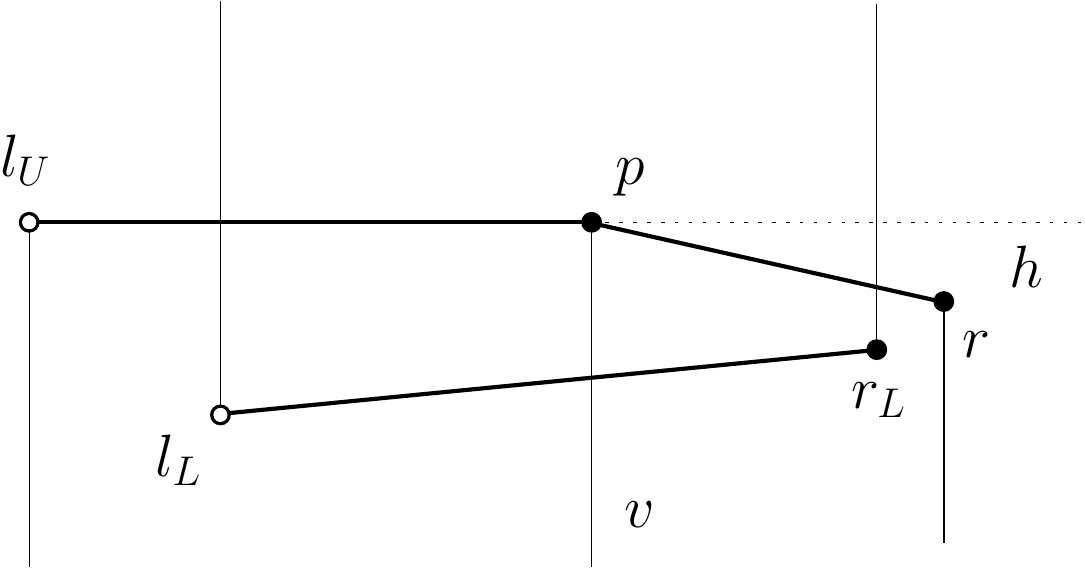}
  }
	\caption{(a) Case (\ref{case:d}) with no point of $\mathcal{P}_U$ to the right of $p$ and no point of $\mathcal{P}_L$ to the left of $l_L$, 
   (b) Case (\ref{case:d})  with the point $r$ of $\mathcal{P}_U$ to the right of $p$}
\end{figure}

% \begin{figure}
% \centering
% \includegraphics[scale=0.5]{FinalA2}
% \caption{ $r_L$ is below the line $pr_U$}
% \label{fig:casea}
% \end{figure}

If none of (\ref{case:a1})-(\ref{case:a4}) occur in case (\ref{case:a}), $r_L$ is below the line $pr_U$. We can assume that the line through $l_L$ and $r_L$ does not intersect the segment $pr_U$ and
that the line through $p$ and $r_U$ does not intersect the segment $l_Lr_L$, see Figure~\ref{fig:casea}. Indeed, otherwise we could reduce this case  (after reversing the $x$ or $y$-axis) to case (\ref{case:b})
with $p$ or $r_L$ playing the role of $p$. Similarly, we can assume that $r_U>r_L$, as otherwise we could use the argument of case (\ref{case:c}).

We color $l_L$ with $blue$ and we color $r_L$ and $r_U$ with $red$.
 We color the points $r\in\mathcal{P}_U$, $r>r_U$ and
 $r\in\mathcal{P}_L$, $r<r_L$, with $blue$.
  We color the remaining points $r\in\mathcal{P}_U$, $r>p$, with $red$.
 Finally, we apply either Observation~\ref{lemma:MainLemma1} or~\ref{lemma:MainLemma2} with $p$ as $p$ and $r_L$ as $q$. Now we finish by arguing
 that our coloring is good.

 A line $l$ witnessing a monochromatic hyperedge cannot avoid both rays corresponding to $p$ and $r_L$.
If $l$ passes below $p$, it also has to pass above $r_U$ and below $r_L$, and we are done by either Observation~\ref{lemma:MainLemma1} or~\ref{lemma:MainLemma2}.
If $l$ passes above $r_L$, it has to pass below $l_L$. Hence, we are again  done by either Observation~\ref{lemma:MainLemma1} or~\ref{lemma:MainLemma2} as well.

\item

% \begin{figure}
% \centering
% \includegraphics[scale=0.5]{FinalD}
% \caption{Case (\ref{case:d})}
% \label{fig:cased}
% \end{figure}

\label{case:d}
In this case we have:  $r_L$ is below $h$, $l_Up$ and $l_Lr_L$ do not intersect each other, and
$l_L>l_U$, as in Figure~\ref{fig:cased0} and~\ref{fig:cased}.

First, we assume that there exist no point of $\mathcal{P}_U$ to the right of $p$ and no point of $\mathcal{P}_L$ to the left of $l_L$, see Figure~\ref{fig:cased0}.
We color the points as follows: $\chi(p)=\chi(r_L)=blue$ and $\chi(l_U)=\chi(l_L)=red$.
We color the points $p'\in\mathcal{P}_U$, $p'<p$, and $p'\in\mathcal{P}_L$, $p'>r_L$, with $red$.
We color the remaining points %$p'\in\mathcal{P}_U$, $p'<l_U$ and
 $p'\in\mathcal{P}_L$, $p'>l_L$, with $blue$.
We can assume that the line through
 $l_Lr_L$ does not intersect the segment $l_Up$, as otherwise we can reduce this case to case (\ref{case:b}) with $l_L$ playing the role of $p$.
 Thus, a line corresponding to a $blue$ monochromatic hyperedge cannot pass above $l_U$. Similarly, a line corresponding to a $red$ monochromatic hyperedge
 cannot pass below $r_L$.

Otherwise, there exists a point $r$, which is either in $\mathcal{P}_U$, such that $r>p$, or in $\mathcal{P}_L$, such that $r<l_L$.
We can assume by symmetry that there is a point $r\in\mathcal{P}_U^0$ immediately to the right of $p$ on $\uhull(\mathcal{P}_U)$, such that $r>r_L$ and  $r$ lies
 below the line through $p$ and $r$, see Figure~\ref{fig:cased}. Indeed, otherwise we could reduce the situation to previous cases (\ref{case:b})-(\ref{case:a}).  We color $r,p$ and $r_L$ with $blue$, and the rest of the points with $red$.
 Our coloring is good since every non-vertical line has to intersect a ray whose corresponding point is $blue$, but no line intersecting all three such rays
 can pass below $l_L$ and above $l_U$.
\end{enumerate}
\noindent
Finally, the case $|\mathcal{P}_U|=1$ can be handled as a special case of
 case (\ref{case:a4}), and that concludes the proof of Theorem~\ref{thm:Main}, except for the algorithmic time complexity.

Most of the proof can be straightforwardly implemented algorithmically.
The bottle neck operations of the algorithm are constructing a convex hull (see e.g \cite{Skiena}), sorting
the points in $\mathcal{P}$ according to the $x$-coordinate and according to their order around a point,
constructing a tangent to a compact convex polygon from a given point, and  the subroutine
from \cite{Keszegh} (in case when there is an uncovered point of the plane).
Since each of these operations takes $O(|V|\log |V|)$ running time, and each of them
is carried out constantly many times, the rest of the theorem follows.

\section{Discussion}
We can generalize our problem as follows.
Let us define $p_{\mathcal{\tilde{H}}}(k)$  as the minimum number $l$ so that
we can $k$-color any finite set of half-planes such that any region covered by at least $l$ half-planes is covered by half-planes of all $k$ colors. In other words, $p_{\mathcal{\tilde{H}}}(k)$ is the minimum number $l$ so that any $l$-fold cover of a set $S\subseteq \mathbb{R}^2$
by a finite set of half-planes can be decomposed into $k$ covers of $S$.

In \cite{Aloupis} it was proved that $p_{\mathcal{\tilde{H}}}(k)\leq 8k-3$, which was recently improved to $4k-3$ in \cite{Smorodinsky2}.
We strengthened this result in one special case, i.e. we proved $p_{\mathcal{\tilde{H}}}(2)=3$.
Thus, it remains as an interesting open question what the right value of $p_{\mathcal{\tilde{H}}}(k)$ is for $k>2$.

\section{Acknowledgment}

\indent The author would like to thank to Bernd G\"{a}rtner, Andreas Razen, and Tibor Szab\'{o} for
 discussions about the problem, and ideas, which improved the paper.
 Thanks are also extended to J\'{a}nos Pach for helpful  suggestions. Last but not least I would
 like to thank David Pritchard for improving the readability of the paper.

%\newpage

%\newpage
%\appendix{{\bf Appendix}}

\bibliographystyle{plain}

\bibliography{23wcf}

\begin{thebibliography}{1}

\bibitem{Aloupis}
Greg Aloupis, Jean Cardinal, S{\'e}bastien Collette, Stefan Langerman, and
  Shakhar Smorodinsky.
\newblock Coloring geometric range spaces.
\newblock {\em Discrete {\&} Computational Geometry}, 41(2):348--362, 2009.

\bibitem{Keszegh}
Bal{\'a}zs Keszegh.
\newblock Weak conflict-free colorings of point sets and simple regions.
\newblock In {\em CCCG}, pages 97--100, 2007.

\bibitem{ll}
L\'{a}szlo Lov\'{a}sz.
\newblock Coverings and colorings of hypergraphs.
\newblock In {\em Proc. 45th Southeastern Conf. Combinatories}, pages 3--12,
  1973.

\bibitem{lll}
L\'{a}szlo Lov\'{a}sz.
\newblock {\em Combinatorial problems and exercises}.
\newblock North-Holland Publishing Company and Akad\'{e}miai Kiad\'{o}, 1992.

\bibitem{Pach}
J\'{a}nos Pach, G\'{a}bor Tardos, and G\'{e}za T\'{o}th.
\newblock Indecomposable coverings.
\newblock In {\em Proceedings of the 7th China-Japan conference on Discrete
  geometry, combinatorics and graph theory}, CJCDGCGT'05, pages 135--148,
  Berlin, Heidelberg, 2007. Springer-Verlag.

\bibitem{Pach2}
János Pach and Gábor Tardos.
\newblock Coloring axis-parallel rectangles.
\newblock {\em Journal of Combinatorial Theory}, 117:776--782, 2010.

\bibitem{Skiena}
Steven~S. Skiena.
\newblock {\em The algorithm design manual}.
\newblock Springer-Verlag New York, Inc., New York, NY, USA, 1998.

\bibitem{Smorodinsky2}
Shakhar Smorodinsky and Yelena Yuditsky.
\newblock Polychromatic coloring for half-planes.
\newblock In {\em SWAT}, pages 118--126, 2010.

\end{thebibliography}

\end{document}